\documentclass[10pt,a4paper]{article}
\usepackage{bbm}
\usepackage{amsfonts}
\usepackage{}
\usepackage{amsmath}
\usepackage{latexsym}
\usepackage{amssymb}
\usepackage{amstext}
\usepackage{fancyhdr}
\usepackage{indentfirst}
\usepackage{graphicx}
\usepackage{epsfig}
\usepackage{pgf}
\usepackage{listings}
\usepackage{cases}
\usepackage{epstopdf}
\usepackage{graphicx}
\usepackage{enumerate}
\usepackage{amsthm}
\usepackage{mathrsfs}
\usepackage{abstract}
\usepackage{caption}
\date{}

\begin{document}
\title{Quantum state transfer on $\mathcal{Q}$-graphs}
\author{Xiao-Qin Zhang$^a$, Shu-Yu Cui$^{b}$, Gui-Xian Tian$^{a, }$\footnote{Corresponding author. E-mail address: gxtian@zjnu.cn or guixiantian@163.com.}\\
    {\small{\it $^a$Department of Mathematics,}}
    {\small{\it Zhejiang Normal University, Jinhua, Zhejiang, 321004, P.R. China}}\\
    {\small{\it $^b$Xingzhi College, Zhejiang Normal University, Jinhua, Zhejiang, 321004, P.R. China}}
}\maketitle
\begin{abstract}
We study the existence of quantum state transfer in $\mathcal{Q}$-graphs in this paper. The $\mathcal{Q}$-graph of a graph $G$, denoted by $\mathcal{Q}(G)$, is the graph derived from $G$ by plugging a new vertex to each edge of $G$ and joining two new vertices which lie on adjacent edges of $G$ by an edge. We show that, if all eigenvalues of a regular graph $G$ are integers, then its $\mathcal{Q}$-graph $\mathcal{Q}(G)$ has no perfect state transfer. In contrast, we also prove that the $\mathcal{Q}$-graph of a regular graph has pretty good state transfer under some mild conditions. Finally, applying the obtained results, we also exhibit many new families of $\mathcal{Q}$-graphs having no perfect state transfer, but admitting pretty good state transfer.


\end{abstract}

\section{Introduction}

The problem of state transfer in quantum walks is a graph-theoretic problem deriving from quantum computing theory. To develop quantum algorithms, Farhi and Gutman in \cite{Farhi} came up with the key concept of continuous-time quantum walk in $1998$. Let $A(G)$ be the adjacency matrix of a graph $G$, then the \emph{transition matrix} of continuous-time quantum walk relative to adjacency matrix $A(G)$ on $G$ is defined by the unitary matrix
\begin{equation}\label{1}
H_{A(G)}(t)=\text{exp}(-itA(G)),
\end{equation}
where $t\in \mathbb{R}$ and $i^2=-1$.
In $2003$, Bose \cite{Bose} studied the problem of information transition in a quantum spin system. Afterwards, Christandl et al. \cite{Christandl} demonstrated that this problem can be reduced to perfect state transfer. Let $e_{u}^n$ denote the characteristic vector of order $n$ corresponding to the vertex $u$ of $G$. Without causing any misunderstanding, $e_{u}^n$ is abbreviated as $e_{u}$. Assume that $u$ and $v$ are two vertices of $G$. If
\begin{equation*}
\text{exp}(-itA(G))e_{u}=\lambda e_{v}
\end{equation*}
with $ |\lambda| =1$, then it is said to admit \textit{perfect state transfer}(PST for short) on $G$ between two vertices $u$ and $v$ at time $t$. This concept has played a crucial role in quantum information and quantum algorithm. Furthermore, if $v=u$, then it is called \emph{periodic} at vertex $u$. Then the topic of characterizing graphs with perfect state transfer has been extensively studied among physic and mathematics communities. In graph theory, it however is well known that graphs admitting PST are rare. A new crucial concept, called pretty good state transfer, was put forward by Godsil in \cite{Godsil2012}, whose restricted condition is more relaxing in comparison to PST. If
\begin{equation*}
|\text{exp}(-itA(G))_{uv}|>1-\epsilon,
\end{equation*}
for any $\epsilon>0$, then it is called as having \emph{pretty good state transfer}(PGST for short) on $G$ between vertices $u$ and $v$ at time $t$. Remark that, if the adjacency matrix is replaced with (signless)Laplacian matrix, then (signless)Laplacian PST and (signless)Laplacian PGST on $G$ can be defined similarly.

Composite graphs are obtained by some operations, such as Cartisian product, the direct product, the strong product and so on, originating from some more simpler graphs.  As composite graphs are more meaningful towards studying larger graphs, it is quite necessary to calculate spectra of compound graphs in many graph-theory problems. Studying state transfer on composite graphs formed by some graph operations tends to be more and more popular in graph-theory field. There are many articles about graphs formed by various graph operations to characterize whether they admit PST (Laplacian PST, signless Laplacian PST ) or PGST (Laplacian PGST, signless Laplacian PGST) or not in recent years. For example, Coutinho and Godsil \cite{Coutinho2016}, using graph products and double covers of graphs, constructed many new graphs admitting PST. Ge et al. \cite{Ge Yang}, using the operation on generalization of the double cones and other graph products, depicted some new constructions of PST graphs, such as lexicographic products, variants of Cartesian product, glued double cones and so on. Angeles-Chaul et al. \cite{Angeles-Canul2010} also obtained many composite graphs admitting PST formed by join operations on graphs and its circulant generalizations. Ackelsberg et al. \cite{Ackelsberg2017} proved that corona product of two simpler graphs under some restrictions has no PST, but some sufficient conditions on corona product admitting PGST were presented. In 2016, Alvir et al. \cite{Alvir2016} investigated Laplacian PST on join operations of graphs, double cones, Cartesian products of graphs and so on. Ackelsberg et al. \cite{Ackelsberg2016} proved that the corona product of two graphs has no Laplacian PST, but it has Laplacian PGST under some special conditions. Li and Liu \cite{Yipeng Li} showed that, the $\mathcal{Q}$-graph of an $r$-regular graph $G$ has no Laplacian PST when $r+1$ is a prime number, but it has Laplacian PGST. Liu and Wang \cite{Xiaogang Liu} proved that the total graph of an $r$-regular graph $G$ exhibits no Laplacian PST when $r+1$ is not a Laplacian eigenvalue of $G$. They also gave a sufficient condition for the total graph admitting Laplacian PGST. In 2021, Tian et al. \cite{Tian2021} demonstrated that there is no signless Laplacian PST on corona product of two graphs, $G\circ K_{m}$, if $m$ equals one or a prime number. They also proved, the corona product $K_{2}\circ \overline{K}_{m}$ has no signless Laplacian PST between two vertices of $K_{2}$ if $m$ is an even number. However, $G\circ \overline{K}_{m}$ exhibits signless Laplacian PGST under some suitable conditions.

Motivated by above some results, we mainly focus on quantum state transfer on $\mathcal {Q}$-graph. The $\mathcal {Q}$-\textit{graph} \cite{Cvetkovic2010} of a graph $G$, denoted by $\mathcal {Q}(G)$, is the graph deriving from $G$ by inverting one new vertex in each edge and joining these new vertices by edges which lie on adjacent edge in $G$. In our work, We show that, if all eigenvalues of a regular graph $G$ are integers, then $\mathcal{Q}(G)$ has no PST. In contrast, we also present a sufficient condition for the $\mathcal {Q}$-graph admitting PGST. Finally, for some distance regular graphs, we also exhibit many new families of their $\mathcal{Q}$-graphs having no PST, but admitting PGST.

\section{Preliminaries}

In whole paper, all graphs considered are finite, simple, undirected and connected. We denote the identity matrix, all-one square matrix of order $n$ and all-one matrix of order $n\times m$ by $I_{n}$, $J_{n}$ and $J_{n,m}$, respectively. For convenience, $j_{n}$ denote the all-one column vector of order $n$. Let $[n]$ be the set of $\{1,2,\ldots, n\}$.

For a graph $G$ of order $n$, let $A(G)$ and $D(G)$ denote its adjacent matrix and degree diagonal matrix, respectively. The eigenvalues of $A(G)$ are called the \textit{eigenvalues} of $G$. Let $\lambda_{0}>\lambda_{1}>\cdots>\lambda_{d}$ be all distinct eigenvalues of $A(G)$ with corresponding multiplicities $a_{0}, a_{1}, \ldots, a_{d}$, where $a_{0}+ a_{1}+\cdots+ a_{d}= |V(G)|= n$. The spectrum of $G$ is denoted by $\text{Sp}(G)$. For each eigenvalue $\lambda_{i}\in \text{Sp}(G)$, let $\{x_{1}^{(i)}, \ldots, x_{a_{i}}^{(i)}\}$ be the set of unit orthogonal eigenvectors of $G$ for $i\in [d]\cup\{0\}$. The transpose of a non-zero column vector $x$ is denoted by $x^T$. It is well-known that the eigenprojector relative to the eigenvalue $\lambda_{i}$ can be written as the following matrix:
\begin{equation*}
{E_{{\lambda_i}}} = \sum\limits_{j = 1}^{{a_i}} {x_j^{(i)}(x_j^{(i)}} {)^T},
\end{equation*}
and $\sum_{i=0}^{d} E_{{\lambda_i}} =I_n$. It is not difficult to see that the $E_{{\lambda_i}} $ is idempotent matrix $(E_{{\lambda_i}}^{2}=E_{{\lambda_i}})$), and $E_{{\lambda_i}}E_{{\lambda_k}}=0$ for $i\neq k$. We can obtain the \textit{spectral decomposition} of $A(G)$ in terms of its eigenprojector:
\begin{equation}\label{2}
A(G)=A(G)\sum_{i=0}^{d} E_{{\lambda_i}}=\sum_{i=1}^{d}\sum_{j=1}^{a_{i}} A(G) x_{j}^{(i)} (x_{j}^{(i)})^T=\sum_{i=0}^{d}\sum_{j=1}^{a_{i}} \lambda_{i} x_{j}^{(i)} (x_{j}^{(i)})^T=\sum_{i=1}^{d}\lambda_{i} E_{{\lambda_i}}.
\end{equation}
On the basis of (\ref{1}) and (\ref{2}), we obtain the alternative expression of the unitary matrix as follows:
\begin{equation}\label{3}
H_{A(G)}(t)=\text{exp}(-itA(G))=\sum_{k\geq 0}\frac{(-i)^{k}A(G)^{k}t^{k}}{k!}=\sum_{i=1}^{d} \text{exp}(-it \lambda_{i})E_{{\lambda_i}}.
\end{equation}

For a vertex $u$ of a graph $G$, we denote its \textit{eigenvalue support} in $G$ by $S=\text{supp}_{A(G)}(u)$, which is the set of all eigenvalues of $A(G)$ satisfying $E_{{\lambda_i}}e_{u}\neq 0$. Two vertices $u$ and $v$ are called as \textit{strongly cospectral} when $E_{\lambda}e_{u}=\pm E_{\lambda}e_{v}$ for every eigenvalue $\lambda$ of $A(G)$. Let $S^{+}$ denote the set of all eigenvalues such that $E_{\lambda}e_{u}=E_{\lambda}e_{v}$, and denote the set of all eigenvalues such that $E_{\lambda}e_{u}=-E_{\lambda}e_{v}$ by $S^{-}$.

The following are some main theorems and lemmas, which will help us to study PST or PGST of $\mathcal {Q}$-graphs.

\paragraph{Theorem 2.1.}(Coutinho \cite{Coutinho2014}) Assume that $G$ is a graph with the vertex set satisfying $|V(G)|\geq 2$, and $u, v\in V(G)$.
If $\lambda_{0}$ is the maximum eigenvalue of $G$,
then $G$ admits PST between the vertices $u$ and $v$ if and only if the following conditions hold.
\begin{enumerate}[(i)]
 \item Two vertices, $u$ and $v$, are strongly cospectral.
 \item Non-zero elements in $\text{supp}_{A(G)}(u)$ are either all integers or all quadratic integers. Moreover, for each eigenvalue $\lambda\in \text{supp}_{A(G)}(u)$, there exists a square-free integer $\Delta$ and integers $a$, $b_\lambda$ such that,
  \begin{equation*}
  \lambda=\frac{1}{2}(a+b_\lambda\sqrt{\Delta}).
  \end{equation*}
  Here we allow $\Delta=1$ if all eigenvalues in $\text{supp}_{A(G)}(u)$ are integers, and $a=0$ if all eigenvalues in $\text{supp}_{A(G)}(u)$ are all multiples of $\sqrt{\Delta}$.
 \item $\lambda\in{S^+}$ if and only if $\dfrac{\lambda_0-\lambda}{g\sqrt{\Delta}}$ is even and  $\lambda\in{S^-}$ if and only if $\dfrac{\lambda_0-\lambda}{g\sqrt{\Delta}}$ is odd, where
 $$g=\gcd\left(\left\{\dfrac{\lambda_0-\lambda}{\sqrt{\Delta}}:\lambda\in \text{supp}_{A(G)}(u)\right\}\right).$$
\end{enumerate}
Moreover, if the above conditions hold, then the following also hold.
\begin{enumerate}[(1)]
     \item There exists a minimum time $\tau_0>0$ at which PST occurs between $u$ and $v$, and
     \begin{equation*}
     \tau_0=\frac{1}{g}\dfrac{\pi}{\sqrt{\Delta}}.
     \end{equation*}
     \item The time of PST $\tau$ is an odd multiple of $\tau_0$.
     \item The phase of PST is given by $\lambda=e^{-it\lambda_0}$.
\end{enumerate}

In order to characterize graphs admitting PST(or PGST), the following two lemmas play a crucial role in our study process.
\paragraph{Lemma 2.2.}(Godsil \cite{Godsil}) If a graph $G$ admits PST between two vertices $u$ and $v$ at time $t$, then $G$ is periodic at vertex $u$ (or $v$) at time $2t$.

\paragraph{Lemma 2.3.}(Godsil \cite{equation}) A graph $G$ at vertex $v$ is periodic if and only if one of the following conditions holds:
\begin{enumerate}[(i)]
\item all element of $\text{supp}_{A(G)}(v)$ are integers;
\item for each eigenvalue of $\text{supp}_{A(G)}(v)$, there is a square-free integer $\Delta$, integer $a$ and corresponding some integer $b_{\lambda}$ so that $\lambda=\frac{1}{2}(a+b_{\lambda}\sqrt{\Delta})$.
\end{enumerate}

\paragraph{Theorem 2.4.}(Hardy and Wright \cite{numbertheory}) Assume that $1, \lambda_1, \ldots, \lambda_m$ are linearly independent over $\mathbb{Q}$. Then, for any real numbers $\alpha_1, \ldots, \alpha_m$ and $N>0$, $\epsilon>0$, there exist integers $\alpha>N$ and $\gamma_1, \ldots, \gamma_m$ such that
\begin{equation}\label{4}
\vert{\alpha\lambda_{k}-\gamma_k-\alpha_k}\vert<\epsilon,
\end{equation}
for each $k\in[m]$. Equivalently, (\ref{4}) can be restated by $\alpha\lambda_{k}-\gamma_k\approx \alpha_k$ for omitting the dependence on $\epsilon$.

\paragraph{Lemma 2.5.}(Richards \cite{galoistheory}) The set $\{\sqrt{\Delta}:\Delta$ is a square-free integer$\}$ is linearly independent over
$\mathbb{Q}$.

\paragraph{Lemma 2.6.}(Coutinho \cite{Coutinho2014}) A real number
$\lambda$ is a quadratic integer if and only if there exist integers $a$, $b$ and $\Delta$ such that $\Delta$ is square-free and one of the following cases holds:
\begin{enumerate}[(i)]
    \item $\lambda=a+b\sqrt{\Delta}$ and $\Delta\equiv2,3\;(\text{mod}\;4)$;
    \item $\lambda=\frac{1}{2}(a+b\sqrt{\Delta}),\; \Delta\equiv1\;(\text{mod}\;4)$, and $a$ and $b$ have the same parity.
\end{enumerate}

\section{Eigenvalues and eigenprojectors of $\mathcal {Q}$-graphs}

In this section, we give the eigenvalues and eigenprojectors of $\mathcal {Q}$-graph of a graph relative to adjacency matrix. These results will play an important role in subsequent studies.

Let $G=(V,E)$ be a graph with vertex set $v_{i}\in V(G)$ and edge set $e_{j}\in E(G)$ with $m\geq n$ for $i\in[n]$, $j\in[m]$. The \textit{vertex-edge incidence matrix} of $G$ is an $n\times m$ matrix $R_{G}=(r_{ij})_{n\times m}$, in which $r_{ij}=1$ if the vertex $v_{i}$ is incident to the edge $e_{j}$, otherwise $r_{ij}=0$. Let $\zeta_{1},\zeta_{2},\ldots,\zeta_{\eta}$ be all unit orthogonal vectors such that $R_G\zeta_i=0$, where $\zeta_i$ is the column vector of order $m$. It is well known \cite{Cvetkovic1995} that $\eta=m-n$ when $G$ is non-bipartite, and $\eta=m-n+1$ when $G$ is bipartite.

Next, we first give the eigenvalues and corresponding eigenvectors of $A(\mathcal {Q}(G))$.

\paragraph{Theorem 3.1.} Let $G$ be an $r$-regular non-bipartite connected graph of order $n$, with $m$ edges and $r\geq2$. Also let $r=\lambda_{0}>\lambda_{1}>\cdots>\lambda_{d}$ be all distinct eigenvalues of $G$ and $x_{1}^{(i)},x_{2}^{(i)}, \ldots, x_{a_{i}}^{(i)}$ be the unit orthogonal eigenvectors corresponding to eigenvalue $\lambda_{i}$ with multiplicity $a_{i}$ for $i\in[d]\cup\{0\}$. Assume that $\zeta_{1}, \zeta_{2}, \ldots, \zeta_{m-n}$ are all unit orthogonal vectors such that $R_{G}\zeta_i=0$. Then the spectrum of $\mathcal{Q}$-graph of $G$ consists precisely of the following:
\begin{enumerate}[(i)]
\item $ \lambda_{i\pm}=\frac{r+\lambda_{i}-2 \pm \sqrt{(\lambda_{i}+r)^{2}+4}}{2}$ are the eigenvalues of $\mathcal {Q}$-graph and the corresponding unite orthogonal eigenvectors are
	\[
	{X_{i \pm }^{j}} = \frac{1}{{\sqrt {{{({\lambda_{i\pm}}+2-r- {\lambda_i})}^2+{\lambda_i+r}}} }}\left( {\begin{array}{*{20}{c}}
		{({\lambda_{i \pm }} + 2 - r - {\lambda_i}){x_j^{(i)}}}\\
		{R_G^T{x_j^{(i)}}}
		\end{array}} \right)
	\]
	for $j\in[a_i]$ and $i\in[d]\cup\{0\}$.
	\item $-2$ is the eigenvalue of $\mathcal{Q}$-graph, with multiplicity $m-n$ and the corresponding orthogonal eigenvectors are
	\[
	Y_k=\frac{1}{||{\zeta_k}||}\left({\begin{array}{*{20}{c}}
	  0\\
	  {{\zeta_k}}
	  \end{array}}\right),
	\]
	for $k\in[m-n]$.

\end{enumerate}

\begin{proof}
In accordance with the definition of $\mathcal {Q}$-graph, then the adjacency matrix of $\mathcal {Q}$-graph of $G$ is given by
\begin{equation*}
{A(\mathcal {Q}(G)}) = \left( {\begin{array}{*{20}{c}}
{O}&{{R_G}}\\
{{R_G}^T}&{A(\ell (G))}
\end{array}} \right),
\end{equation*}
where $\ell(G)$ denotes the line graph of $G$. Remark that the characteristic polynomial of $A(\mathcal {Q}(G))$ relative to the characteristic polynomial of $A(\ell(G))$ has been obtained in \cite{Cvetkovic1975}. Here we give the detailed proof for the convenience of readers. Since $A(\ell(G))=R_{G}^{T}R_{G}-2I_{m}$ and $Q_{G}=R_{G}R_{G}^{T}$, where $Q_{G}$ is the signless Laplacian matrix of $G$. Then the characteristic polynomial of $\mathcal {Q}(G)$
\begin{equation*}
\begin{split}
{P_{{A(\mathcal {Q}(G))}}}(t) &=\det \left( {\begin{array}{*{20}{c}}
{t{I_n}}&{ - {R_G}}\\
{ - R_G^T}&{t{I_m} - A(\ell (G))}
\end{array}} \right)\\
& = \det\left( {\begin{array}{*{20}{c}}
{t{I_n}}&{ - {R_G}}\\
{ - R_G^T}&{(t + 2){I_m} - R_G^T{R_G}}
\end{array}} \right)\\
& =\det \left( {\begin{array}{*{20}{c}}
{t{I_n}}&{ - {R_G}}\\
{(- t - 1)R_G^T}&{(t + 2){I_m}}
\end{array}} \right)\\
&=\det \left( {\begin{array}{*{20}{c}}
{t{I_n} - \frac{{ t + 1}}{{t+ 2}}{R_G}R_G^T}&0\\
{(-t - 1)R_G^T}&{(t + 2 ){I_m}}
\end{array}} \right)\\
& = {(t + 2 )^m}\det (t{I_n} - \frac{{ t + 1}}{{t + 2}}{R_G}R_G^T)\\
& = {(t + 2 )^{m - n}}\det ((t + 2)t{I_n} - (t + 1){R_G}{R_G}^T)\\
&= {(t + 2 )^{m - n}}\det ((t + 2)t{I_n} - (t + 1){Q_G}).
\end{split}
\end{equation*}
Since $Q_G=rI_n+A(G)$, then the characteristic polynomial
\[
{P_{{A(\mathcal {Q}(G))}}}(t) = {(t + 2 )^{m - n}} \prod \limits_{i = 0}^d {[(t + 2 )t - (t + 1){(\lambda_i+r)}]^{{a_i}}},
\]
which implies that $-2$ and
\[
{\lambda_{i \pm }} = \frac{{r + {\lambda_i} - 2 \pm \sqrt {{{(\lambda_{i} + r )}^2} + 4} }}{2}
\]
are the eigenvalues of $\mathcal {Q}(G)$ for $i\in[d]\cup\{0\}$. Since $G$ is a non-bipartite graph, then there is no $\lambda_{i}$ such that $\lambda_i= -r$ for any $i\in[d]$. This implies that $\lambda_{i \pm }\neq-2$ for $i\in[d]\cup\{0\}$. Hence, $-2$ and $\lambda_{i \pm }$ are the eigenvalues of $A(\mathcal {Q}(G))$ with corresponding respective multiplicities $m-n$ and $a_{i}$ for $i\in[d]\cup\{0\}$.

Next, we give all eigenvectors of $A(\mathcal {Q}(G))$. First, by a simple calculation, we have
\[
A(\mathcal {Q}(G))Y_k=
\left( {\begin{array}{*{20}{c}}
{O}&{{R_G}}\\
{{R_G}^T}&{ A(\ell (G))}
\end{array}} \right)\frac{1}{||{\zeta_k}||}\left( {\begin{array}{*{20}{c}}
0\\
\zeta_k
\end{array}} \right) = - 2\frac{1}{||{\zeta_k}||}\left( {\begin{array}{*{20}{c}}
0\\
\zeta_k
\end{array}} \right)= - 2Y_k,
\]
which implies that, for every $ \zeta_k \in \{ \zeta_1, \zeta_2, \ldots, \zeta_{m - n}\}$,
\[
Y_k = \frac{1}{||{\zeta_k}||}\left( {\begin{array}{*{20}{c}}
0\\
\zeta_k
\end{array}} \right)
\]
is the eigenvector corresponding to eigenvalue $-2$. Since $\{x_1^{(i)}, x_2^{(i)}, \ldots, x_{{a_i}}^{(i)}\}$ is an unit orthogonal basis of the eigenspace $V_{\lambda_{i}}$ corresponding to the eigenvalue $\lambda_{i}$. Then $A(G)x_{j}^{(i)}=\lambda_{i}x_{j}^{(i)}$ for $j=[a_{i}]$ and $i= [d]\cup\{0\}$. It is easy to see that
\[
{\lambda_{i \pm }}= \frac{{{\lambda_i+r}}}{{{\lambda_{i \pm }} + 2 - r - {\lambda_i}}}.
\]
Thus we obtain
\begin{equation*}\small
\begin{split}
{A(\mathcal {Q}(G))}X_{i \pm }^j&=\frac{1}{{\sqrt {{{({\lambda_{i \pm }} + 2 - r - {\lambda_i})}^2} + {\lambda_i}+r} }}\left( {\begin{array}{*{20}{c}}
{O}&{{R_G}}\\
{{R_G}^T}&{ A(\ell (G))}
\end{array}} \right)\left( {\begin{array}{*{20}{c}}
{({\lambda_{i \pm }} + 2 - r - {\lambda_i})x_j^{(i)}}\\
{R_G^Tx_j^{(i)}}
\end{array}} \right)\\
&= \frac{1}{{\sqrt {{{({\lambda_{i \pm }} + 2 - r - {\lambda_i})}^2} + {\lambda_i}+r} }}\left( {\begin{array}{*{20}{c}}
{O}&{{R_G}}\\
{{R_G}^T}&{- 2{I_m} + R_G^T{R_G}}
\end{array}} \right)\left( {\begin{array}{*{20}{c}}
{({\lambda_{i \pm }} + 2 - r - {\lambda_i})x_j^{(i)}}\\
{R_G^Tx_j^{(i)}}
\end{array}} \right)\\
&= \frac{1}{{\sqrt {{{({\lambda_{i \pm }} + 2 - r - {\lambda_i})}^2} + {\lambda_i}+r} }}\left( {\begin{array}{*{20}{c}}
{(\lambda_{i}+r)x_j^{(i)}}\\
{({\lambda_{i \pm }} + 2 - r - {\lambda_i} + r - 2 + {\lambda_i})R_G^Tx_j^{(i)}}
\end{array}} \right)\\
&= \frac{1}{{\sqrt {{{({\lambda_{i \pm }} + 2 - r - {\lambda_i})}^2} + {\lambda_i+r}} }}{\lambda_{i \pm }}\left( {\begin{array}{*{20}{c}}
{({\lambda_{i \pm }} + 2 - r - {\lambda_i})x_j^{(i)}}\\
{R_G^Tx_j^{(i)}}
\end{array}} \right)\\
&= {\lambda_{i \pm }}X_{i \pm }^j.
\end{split}
\end{equation*}
Hence, $X_{i\pm}^{j}$ are the eigenvectors of $A(\mathcal {Q}(G))$ corresponding to the  eigenvalues $\lambda_{i\pm}$.

In what follows, we still need to prove that all $X_{i\pm}^{j}$ and $Y_k$ are orthogonal.
From the process of proof above, $(X_{i\pm}^{j})^TX_{i\pm}^{k}=0$ and $Y_j^{T}Y_k=0$ for $j\neq k$. Clearly, $(X_{i\pm}^{j})^TY_{k}=0$ for $j\in [a_{i}]$, $i= [d]\cup\{0\}$ and $k\in[m-n]$.
Thus, we will only prove that $(X_{i+ }^j)^TX_{i- }^j=0$. As a matter of convenience, set
$$\alpha=\frac{1}{{\sqrt {{{({\lambda_{i + }} + 2 - r - {\lambda_i})}^2+{\lambda_i+r}}} }}\frac{1}{{\sqrt {{{({\lambda_{i- }} + 2 - r - {\lambda_i})}^2+{\lambda_i+r}}} }}.$$
Since $({\lambda_{i + }} + 2 - r - {\lambda_i})({\lambda_{i - }} + 2 - r - {\lambda_i}) =  - (\lambda_i+r)$, then
\begin{equation*}
\begin{split}
{(X_{i + }^j)^T}X_{i - }^j &= \alpha{\left( {\begin{array}{*{20}{c}}
{({\lambda_{i _{+} }} + 2 - r - {\lambda_i})x_j^{(i)}}\\
{R_G^Tx_j^{(i)}}
\end{array}} \right)^T}\left( {\begin{array}{*{20}{c}}
{({\lambda_{i _{-} }} + 2 - r - {\lambda_i})x_j^{(i)}}\\
{R_G^Tx_j^{(i)}}
\end{array}} \right)\\
&= \alpha(({\lambda_{i + }} + 2 - r - {\lambda_i})({\lambda_{i - }} + 2 - r - {\lambda_i}){(x_j^{(i)})^T}x_j^{(i)} + {(x_j^{(i)})^T}{R_G}R_G^Tx_j^{(i)})\\
&=  \alpha(- ({\lambda_i+r}){(x_j^{(i)})^T}x_j^{(i)} + ({\lambda_i+r}){(x_j^{(i)})^T}x_j^{(i)})\\
&= 0.
\end{split}
\end{equation*}
Then $X_{i+}^j$ and $X_{i - }^j$ are orthogonal eigenvectors for any $j\in [a_{i}]$ and $i\in [d]\cup\{0\}$.
Above all, all $X_{i\pm}^j$ and $Y_k$ are orthogonal eigenvectors for $j\in [a_{i}]$, $i= [d]\cup\{0\}$ and $k\in[m-n]$. This completes our proof.
\end{proof}

\paragraph{Theorem 3.2.}
Let $G$ be an $r$-regular $(r\geq2)$ bipartite connected graph of order $n$ with $m$ edges, and $V_{1}\cup V_{2}$ be the bipartition of the vertex set of $G$. Also let $ r=\lambda_{0} > \lambda_{1}> \cdots > \lambda_{d}=-r$ are all different eigenvalues of $G$ with the corresponding multiplicities $a_{0}, a_{1}, \ldots a_{d}$, and $\{x_{1}^{(i)},x_{2}^{(i)}, \ldots, x_{a_{i}}^{(i)}\}$ be the unit orthogonal eigenvectors corresponding to $\lambda_{i}$ for $i\in [d]\cup\{0\}$. Suppose that $\zeta_{1}, \zeta_{2}, \ldots, \zeta_{m-n+1}$ are all unit orthogonal vectors such that $R_{G}\zeta_k=0$ for $k\in[m-n+1]$. Then $-2$, $0$ and
$$ \lambda_{i\pm}=\frac{r+\lambda_{i}-2 \pm \sqrt{(\lambda_{i}+r)^{2}+4}}{2}$$
are the eigenvalues of $\mathcal {Q}$-graph and the corresponding eigenvectors are, for $k=[m-n+1]$,
\[
Y_k = \frac{1}{||{\zeta_k}||}\left( {\begin{array}{*{20}{c}}
	0\\
	{{\zeta_k}}
	\end{array}} \right),\;\;
\frac{1}{{\sqrt n }}\left( {\begin{array}{*{20}{c}}
	{{j_{{\frac{n}{2}}}}}\\
	{ - {j_{\frac{n}{2}}}}\\
	{{0_m}}
	\end{array}} \right)
\]
and
\[
{X_{i \pm }^{j}} = \frac{1}{{\sqrt {{{({\lambda_{i \pm }} + 2 - r - {\lambda_i})}^2+{\lambda_i+r}}} }}\left( {\begin{array}{*{20}{c}}
{({\lambda_{i \pm }} + 2 - r - {\lambda_i}){x_j^{(i)}}}\\
{R_G^T{x_j^{(i)}}}
\end{array}} \right)
\]
for $j\in [a_{i}]$, $i\in [d-1]\cup\{0\}$.

\begin{proof}
Since $G$ is bipartite connected $r$-regular graph, then the smallest eigenvalue $\lambda_{d}=-r$,  $|V_{1}|=|V_{2}|=\frac{n}{2}$ and the multiplicities of $\lambda_{0}$ and $\lambda_{d}$ are equal to 1. In the light of
\[
{\lambda_{i \pm }} = \frac{{r + {\lambda_i} - 2 \pm \sqrt {{{(\lambda_{i} + r)}^2} + 4} }}{2},
\]
we have ${\lambda_{d + }}=0$ and ${\lambda_{d - }} = -2$.
It is not difficult to verify that
\[
{A(\mathcal {Q}(G))}\frac{1}{{\sqrt n }}\left( {\begin{array}{*{20}{c}}
{{j_{\frac{n}{2}}}}\\
{ - {j_{\frac{n}{2}}}}\\
{{0_m}}
\end{array}} \right)\\
 = \left( {\begin{array}{*{20}{c}}
{O}&{{R_G}}\\
{{R_G}^T}&{ - 2{I_m} + R_G^T{R_G}}
\end{array}} \right)\frac{1}{{\sqrt n }}\left( {\begin{array}{*{20}{c}}
{{j_{\frac{n}{2}}}}\\
{ - {j_{\frac{n}{2}}}}\\
{{0_m}}
\end{array}} \right)\\
 = 0\frac{1}{{\sqrt n }}\left( {\begin{array}{*{20}{c}}
{{j_{\frac{n}{2}}}}\\
{ - {j_{\frac{n}{2}}}}\\
{{0_m}}
\end{array}} \right).
\]
Thus
\[\frac{1}{{\sqrt n }}\left( {\begin{array}{*{20}{c}}
{{j_{\frac{n}{2}}}}\\
{ - {j_{\frac{n}{2}}}}\\
{{0_m}}
\end{array}} \right)\]
is the  eigenvector corresponding to the eigenvalue $0$ of $A(\mathcal {Q}(G))$. The rest of the proof is similar exactly to that of Theorem 3.1, omitted.
\end{proof}

On the grounds of Theorems 3.1 and 3.2, we obtain immediately the eigenprojectors and spectral decomposition of $A(\mathcal {Q}(G))$.

\paragraph{Theorem 3.3.} Assume that $G$ is an $r$-regular $(r\geq2)$ connected graph of order $n$ with $m$ edges. Then
\begin{enumerate}[(a)]
\item If $G$ is non-bipartite, then the  eigenprojectors corresponding to the respective eigenvalues $\lambda_{i\pm}$ and $-2$ of $A(\mathcal {Q}(G))$ are denoted by $F_{\lambda_{i\pm}}$ and $F_{-2}$, where
\begin{equation}\small\label{5}
{F_{{\lambda_{i \pm }}}} = \frac{1}{{{{({\lambda_{i \pm }} + 2 - r - {\lambda_i})}^2} + {\lambda_i}+r}}\left( {\begin{array}{*{20}{c}}
{{{({\lambda_{i \pm }} + 2 - r - {\lambda_i})}^2}{E_{{\lambda_i}}}}&{({\lambda_{i \pm }} + 2 - r - {\lambda_i}){E_{{\lambda_i}}}{R_G}}\\
{({\lambda_{i \pm }} + 2 - r - {\lambda_i}){{({E_{{\lambda_i}}}{R_G})}^T}}&{R_G^T{E_{{\lambda_i}}}{R_G}}
\end{array}} \right)
\end{equation}
and
\begin{equation}\label{6}
{F_{- 2}} = \sum\limits_{k = 1}^{m - n} {\frac{1}{{||{\zeta_k}||^{2}}}} \left( {\begin{array}{*{20}{c}}
0&0\\
0&{{\zeta_k}\zeta_k^T}
\end{array}} \right).
\end{equation}
Thus, we get the spectral decomposition of $A(\mathcal {Q}(G))$ as follows:
\begin{equation}\label{7}
{A(\mathcal {Q}(G))} = \sum\limits_{i = 0}^{d} {\sum\limits_ \pm  {{\lambda_{i \pm }}{F_{{\lambda_{i \pm }}}} + (- 2)} } {F_{- 2}}.
\end{equation}
\item If $G$ is bipartite, then the eigenprojectors corresponding to the respective eigenvalues $-2$, $0$ and $\lambda_{i\pm}$ of $A(\mathcal {Q}(G))$ are denoted by $F_{-2}$, $F_{0}$ and $F_{\lambda_{i\pm}}$, where
\begin{equation}\label{8}
{F_{- 2}} = \sum\limits_{k = 1}^{m - n+1} {\frac{1}{{||{\zeta_k}||^{2}}}} \left( {\begin{array}{*{20}{c}}
0&0\\
0&{{\zeta_k}\zeta_k^T}
\end{array}} \right),
\end{equation}
\begin{equation}\label{9}
{F_0} = \frac{1}{n}\left( {\begin{array}{*{20}{c}}
{{J_{\frac{n}{2}}}}&{ - {J_{\frac{n}{2}}}}&0\\
{ - {J_{\frac{n}{2}}}}&{{J_{\frac{n}{2}}}}&0\\
0&0&0
\end{array}} \right) = \left( {\begin{array}{*{20}{c}}
{{E_{-r}}}&0\\
0&0
\end{array}} \right),
\end{equation}
and for $i\neq d$
\begin{equation}\small\label{10}
{F_{{\lambda_{i \pm }}}} = \frac{1}{{{{({\lambda_{i \pm }} + 2 - r - {\lambda_i})}^2} + {\lambda_i+r}}}\left( {\begin{array}{*{20}{c}}
	{{{({\lambda_{i \pm }} + 2 - r - {\lambda_i})}^2}{E_{{\lambda_i}}}}&{({\lambda_{i \pm }} + 2 - r - {\lambda_i}){E_{{\lambda_i}}}{R_G}}\\
	{({\lambda_{i \pm }} + 2 - r - {\lambda_i}){{({E_{{\lambda_i}}}{R_G})}^T}}&{R_G^T{E_{{\lambda_i}}}{R_G}}
	\end{array}} \right).
\end{equation}
Thus, we get the spectral decomposition of $A(\mathcal {Q}(G))$ as follows:
\begin{equation}\label{11}
{A_{\mathcal {Q}(G)}} =\sum\limits_{i = 0}^{d - 1} {\sum\limits_ \pm  {{\lambda_{i \pm }}{F_{{\lambda_{i \pm }}}} + ( - 2)} } {F_{- 2}} + 0{F_0}.
\end{equation}
\end{enumerate}

In accordance with Theorems 3.1, 3.2 and 3.3, we obtain easily the following proposition, which will be applied to analyse quantum state transfer of $\mathcal {Q}(G)$.

\paragraph{Proposition 3.4.} Let $G$ be an $r$-regular $(r\geq2)$ connected graph of order $n\geq2$ with $m$ edges. Then, for two vertices $u,v\in V(G)$, we have
\begin{enumerate}[(i)]
\item If $G$ is non-bipartite, then
\begin{equation}\small\label{12}
{(e_v^{m + n})^T}\text{exp}( - it{A(\mathcal {Q}(G))})e_u^{m + n}
 = {e^{ - it\frac{{r - 2}}{2}}}\sum\limits_{i = 0}^{d} {{e^{ - it\frac{{{\lambda_i}}}{2}}}} e_v^T{E_{{\lambda_i}}}{e_u}(\cos\frac{{{\Delta _{{\lambda_i}}}t}}{2} + \frac{{{\lambda_i} + r - 2}}{{{\Delta _{{\lambda_i}}}}}i\sin\frac{{{\Delta _{{\lambda_i}}}t}}{2}),
\end{equation}
where $\Delta_{\lambda_i}=\sqrt{(\lambda_i+r)^2+4}$ for $i\in [d]\cup \{0\}$.
\item If $G$ is bipartite, then
\begin{equation}\small\label{13}
\begin{split}
(e_v^{m + n})&^T\text{exp}( - it{A(\mathcal {Q}(G))})e_u^{m + n}\\
 &= {e^{ - it\frac{{r - 2}}{2}}}\sum\limits_{i = 0}^{d - 1} {{e^{ - it\frac{{{\lambda_i}}}{2}}}} e_v^T{E_{{\lambda_i}}}{e_u}(\cos\frac{{{\Delta _{{\lambda_i}}}t}}{2} + \frac{{{\lambda_i} + r - 2}}{{{\Delta _{{\lambda_i}}}}}i\sin\frac{{{\Delta _{{\lambda_i}}}t}}{2}) + {e^{ - it0}}e_v^T{E_{-r}}{e_u}.
\end{split}
\end{equation}
where $\Delta_{\lambda_i}=\sqrt{(\lambda_i+r)^2+4}$ for $i\in [d-1]\cup \{0\}$
\end{enumerate}

\begin{proof} (i) Since $G$ is a non-bipartite graph. Then the (a) of Theorem 3.3 and (\ref{3}) imply that
\[
\text{exp}( - it{A(\mathcal {Q}(G))}) = \sum\limits_{i = 0}^{d} {\sum\limits_ \pm  {{e^{ - it{\lambda_{i \pm }}}}{F_{{\lambda_{i \pm }}}} + {e^{ - it(- 2)}}} } {F_{- 2}}.
\]
By a simple calculation, we get
\begin{equation}\label{14}
({\lambda_{i + }} + 2 - r - {\lambda_i})({\lambda_{i - }} + 2 - r - {\lambda_i}) =  - ({\lambda_i+r}),
\end{equation}
\begin{equation}\label{15}
{({\lambda_{i + }} + 2 - r - {\lambda_i})^2} + {({\lambda_{i - }} + 2 - r - {\lambda_i})^2} = {\Delta _{{\lambda_i}}}^2 - 2({\lambda_i+r})
\end{equation}
and
\begin{equation}\label{16}
{({\lambda_{i - }} + 2 - r - {\lambda_i})^2} - {({\lambda_{i + }} + 2 - r - {\lambda_i})^2} = ({\lambda_i} + r - 2){\Delta _{{\lambda_i}}}.
\end{equation}
It follows from the formulas (\ref{5}), (\ref{6}) and (\ref{7}) that, for two vertices $u, v\in V(G)$,
\begin{equation*}
\begin{split}
{(e_v^{m + n})^T}\text{exp}( - it{A(\mathcal {Q}(G))})e_u^{m + n}
 &= {(e_v^{m + n})^T}(\sum\limits_{i = 0}^{d} {\sum\limits_ \pm  {{e^{ - it{\lambda_{i \pm }}}}{F_{{\lambda_{i \pm }}}} + {e^{ - it(- 2)}}} } {F_{- 2}})e_u^{m + n}\\
 &= \sum\limits_{i = 0}^{d} {\sum\limits_ \pm  {{e^{ - it{\lambda_{i \pm }}}}{{(e_v^{m + n})}^T}{F_{{\lambda_{i \pm }}}}} } e_u^{m + n}\\
 &= {e^{ - it\frac{{r - 2}}{2}}}\sum\limits_{i = 0}^{d} {{e^{ - it\frac{{{\lambda_i}}}{2}}}} e_v^T{E_{{\lambda_i}}}{e_u}(\cos\frac{{{\Delta _{{\lambda_i}}}t}}{2} + \frac{{{\lambda_i} + r - 2}}{{{\Delta _{{\lambda_i}}}}}i\sin\frac{{{\Delta _{{\lambda_i}}}t}}{2}).
\end{split}
\end{equation*}

(ii) Assume that $G$ is a bipartite graph. It follows from the (b) of Theorem 3.3 and (\ref{3}) that
\[
\text{exp}( - it{A_{\mathcal {Q}(G)}}) = \sum\limits_{i = 0}^{d} {\sum\limits_ \pm  {{e^{ - it{\lambda_{i \pm }}}}{F_{{\lambda_{i \pm }}}} + {e^{ - it(- 2)}}} } {F_{- 2}} + {e^{ - it0}}{F_0}.
\]
The rest proof is exactly similar to that of the (i), the details are omitted.
\end{proof}

\section{PST on $\mathcal {Q}$-graphs}

We mainly investigate the existence of PST of $\mathcal {Q}$-graph in this section. First we exhibit a lemma, which is given by Li and Liu in \cite{Yipeng Li}.

\paragraph{Lemma 4.1.}(\cite{Yipeng Li}) Let $G$ be an $r$-regular connected graph of order $n$ with $m$ edges and $r\geq2$. Also let all distinct eigenvalues of $A(G)$ be as described at the beginning of Section 2, along with the eigenprojectors $E_{\lambda_{0}}, E_{\lambda_{1}},\ldots,E_{\lambda_{d}}$, respectively.
\begin{enumerate}[(i)]
\item If $G$ is non-bipartite, then there is some index $i_0\in [d]$ such that $(E_{\lambda_{i_0}}R_{G})e_{k}^{m}\neq0$ for any $k\in [m]$.
\item If $G$ is bipartite, then for any $k\in [m]$, there is some index $i_0\in[d-1]$ such that $(E_{\lambda_{i_0}}R_{G})e_{k}^{m}\neq0$.
\end{enumerate}

\begin{proof}
Since $G$ is $r$-regular graph, then $A(G)=D(G)-L(G)=rI_{n}-L(G)$ where $L(G)$ is the Laplacian matrix of $G$. Therefore, all the eigenvectors of $A(G)$ are also the eigenvectors of $L(G)$, which implies that the eigenprojector of $A(G)$ corresponding to the eigenvalue $\lambda_i$ is also the eigenprojector of $L(G)$ corresponding to the eigenvalue $r-\lambda_i$. Hence, we obtained the desired results by the Lemma 4.1 of \cite{Yipeng Li}.
\end{proof}

Lemma 2.2 indicated that, if a graph $G$ admits PST between two vertices $u$ and $w$, then $G$ must be periodic at  two vertices $u$ and $w$ of $G$. By analyzing periodic conditions on $\mathcal {Q}(G)$, we obtain the following theorem.

\paragraph{Theorem 4.2.} Let $G$ be an $r$-regular($r\geq2$) connected graph of order $n>2$ with $m$ edges. Let $V(G)\cup N(G)$ be the vertex set of $\mathcal {Q}(G)$, where $N(G)$ is the new vertex set in $\mathcal {Q}(G)$.
If all eigenvalues of $A(G)$ are integers, then there is no PST in $\mathcal {Q}(G)$.

\begin{proof} Here we only proof the case of non-bipartite graph $G$. For bipartite graph, the proof is similar to that of non-bipartite graph.

Let $z$ be any vertex of $\mathcal {Q}(G)$. Suppose towards the contradiction that $\mathcal {Q}(G)$ has PST between $z$ and another vertex. Lemma 2.2 implies that $\mathcal {Q}(G)$ is periodic at vertex $z$. In the light of Lemma 2.3, we only need to discuss it in the following two cases.

\emph{Case 1.} All the elements of $\text{supp}_{A(\mathcal {Q}(G))}(z)$ are integers, where $z\in V(G)\cup N(G)$.

We first suppose that the vertex $z\in V(G)$. Since $G$ is a connected graph, then there is a eigenvalue $\lambda\in \text{supp}_{A(G)}(z)$. Thus $E_{\lambda}e_{z}\neq 0$ by the definition of eigenvalue support. According to the (a) of Theorem 3.3, $F_{\lambda\pm}e_{z}^{m+n}\neq 0$ as the coefficient ${\lambda_{ \pm }} + 2 - r - {\lambda}$ is non-zero for any $\lambda$ when $G$ is non-bipartite graph. This implies that $\lambda_{\pm}\in\text{supp}_{A(\mathcal {Q}(G))}(z)$.
Since
\[
{\lambda_ + } = \frac{{r + \lambda - 2 + \sqrt {{{(r + \lambda )}^2} + 4} }}{2},\;\;
{\lambda_ - } = \frac{{r + \lambda - 2 - \sqrt {{{(r + \lambda )}^2} + 4} }}{2},
\]
then $\lambda_{+}+\lambda_{-}=r+\lambda-2$, $\lambda_{+}-\lambda_{-}=\sqrt{(r+\lambda)^{2}+4}$. Bear in mind that all of $\lambda$, $\lambda_{+}$ and $\lambda_{-}$ are integers. It follows that $\sqrt{(r+\lambda)^{2}+4}$ is integer, in other words, $(r+\lambda)^{2}+4$ is a perfect square. Observing that  $4$ is even and
\[
{(\lambda + r )^2} < {(\lambda + r )^2} + 4 \leq{(\lambda + r +2)^2}={(\lambda + r )^2}+4+4(\lambda+r).
\]
Since $\sqrt{(r+\lambda)^{2}+4}$ and $\lambda + r $ have the same parity, hence $\sqrt{(r+\lambda)^{2}+4}$ can not be an integer whenever $\lambda+r\neq 0$. This contradicts the fact that all the elements of $\text{supp}_{A(\mathcal {Q}(G))}(z)$ are integers.

Next suppose that the vertex $z\in N(G)$. From the (i) of Lemma 4.1, there exists some eigenvalue $\lambda_{i_{0}}$ such that $(E_{\lambda_{i_{0}}}R_{G})e_{z}^{m}\neq0$ for $i_{0}\in [d]$. Thus, according to the (a) of Theorem 3.3, we obtain that
\[
{F_{{\lambda_{i_{0}\pm} }}}e_z^{m + n} = \frac{1}{{{{({\lambda_{i_{0}\pm} } + 2 - r - \lambda_{i_{0}})}^2} + \lambda_{i_{0}}+r}}\left( {\begin{array}{*{20}{c}}
{({\lambda_{i_{0}\pm} } + 2 - r - \lambda_{i_{0}})({E_{\lambda_{i_{0}}}}{R_G})e_z^m}\\
{(R_G^T{E_{\lambda_{i_{0}}}}{R_G})e_z^m}
\end{array}} \right)\neq0.
\]
So $\lambda_{i_{0}\pm}\in \text{supp}_{A(\mathcal {Q}(G))}(z)$ and $\lambda_{i_{0}\pm}$ are integers. Similar to the proof of the case of $z\in V(G)$ above, it results in that not all the elements of $\text{supp}_{A(\mathcal {Q}(G))}(z)$ are integers, a contradiction.

\emph{Case 2.} All the elements of $\text{supp}_{A(\mathcal {Q}(G))}(z)$ are the form of $\lambda_{\pm}=\frac{a+b\sqrt{\Delta}}{2}$ for integers $a$ and square-free integer $\Delta$ with some integer $b_{\pm}$.

First assume that the vertex $z \in V(G)$. On the basis of the proof of Case $1$ above, we have $\lambda_\pm\in \text{supp}_{A(\mathcal {Q}(G))}(z)$. Since all the elements of $\text{supp}_{A(\mathcal {Q}(G))}(z)$ are the form of $\frac{a+b_{\pm}\sqrt{\Delta}}{2}$ and $(\lambda + r )^2 + 4 $ can not be perfect square for integer $\lambda$. Then, there exists a square-free integer $\Delta>1$, along with integers $a$ (we allow that $a=0$), $b_{+}$ and $b_{-}$, such that
\[
{\lambda_\pm } = \frac{{a + {b_ \pm }\sqrt \Delta  }}{2}.
\]
Thus, it follows from (\ref{14}) that
\begin{equation*}
\frac{1}{4}{b_ + }{b_ - }\Delta  + \frac{1}{4}{(a + 4 - 2r - 2\lambda)^2} + \frac{1}{4}\sqrt \Delta  ({b_ + } + {b_ - })(a + 4 - 2r - 2\lambda) =  - (\lambda+r).
\end{equation*}
Since  $-(\lambda+r)$ is integer and $\sqrt{\Delta}$ is irrational, then either $b_ +  + b_ - =0$ or $a+4-2r-2\lambda=0$. If $b_ +  + b_ - =0$, then $\lambda_{+}+\lambda_{-}=a=r+\lambda-2$. At this time, $|\text{supp}_{A(G)}(z)|=1$, which contradicts that $G$ is a connected graph with $n\geq2$ vertices. Otherwise, $a+4-2r-2\lambda=0$, then $a=2r+2\lambda-4$, which also implies that $|\text{supp}_{A(G)}(z)|=1$, a contradiction. Hence, not all the elements of $\text{supp}_{A(\mathcal {Q}(G))}(z)$ are quadratic integer for $b_+ + b_- \neq0$ and $a+4-2r-2\lambda\neq0$.

Next let the vertex $z\in N(G)$. In the light of the (i) of Lemma 4.1 again, there exists some eigenvalue $\lambda_{i_{0}}$ such that $(E_{\lambda_{i_{0}}}R_{G})e_{z}^{m} \neq 0$. Then, from Theorem 3.3, we obtain that
\[
{F_{{\lambda_{i_{0}\pm} }}}e_z^{m + n} = \frac{1}{{{{({\lambda_ {i_{0}\pm} } + 2 - r - \lambda_{i_{0}})}^2} + \lambda_{i_{0}}+r}}\left( {\begin{array}{*{20}{c}}
{({\lambda_ {i_{0}\pm} } + 2 - r - \lambda_{i_{0}})({E_{\lambda_{i_{0}}}}{R_G})e_z^m}\\
{(R_G^T{E_{\lambda_{i_{0}}}}{R_G})e_z^m}
\end{array}} \right)\neq0.
\]
So, $\lambda_{i_{0}\pm} \in \text{supp}_{A(\mathcal {Q}(G))}(z)$ and $\lambda_{i_{0}\pm}$ are quadratic integers. Similar to the discussion in Case 1, not all the elements of $\text{supp}_{A(\mathcal {Q}(G))}(z)$ are quadratic integers.

Neither of two cases above is likely to happen. Hence, it follows from Lemma 2.3 that there is no PST in $\mathcal {Q}(G)$. This completes the proof of theorem.
\end{proof}

Remark that the condition $r\geq2$ must be required in Theorem 4.2. In fact, if we allow that $r=1$, then $G$ is the path $P_2$ with $\text{Sp}(P_2)=\{-1,1\}$. However, Christandl et al in \cite{Christandl} proved that the $\mathcal {Q}$-graph $P_{3}$ of the path $P_{2}$ has PST at antipodal vertices. In addition, in accordance with the proof of Theorem 4.2, we easily arrive at the following corollary.

\paragraph{Corollary 4.3.} Let $G$ be an $r$-regular connected graph with $r\geq2$. Assume that $u$ is any vertex in $G$.
If all the elements of $\text{supp}_{A(G)}(u)$ are integers, then $\mathcal {Q}(G)$ has no PST between $u$ and any other vertex of $\mathcal {Q}(G)$.

\paragraph{Example 1.} Assume that $G$ is one of the following graph:
\begin{enumerate}[(i)]
\item $d$-dimensional hypercubes $Q_{d}$ for $d\geq2$;

\item Cocktail party graphs $\overline{mK_{2}}$ for $m\geq 2$;

\item Halved $2d$-dimensional hypercubes $\frac{1}{2}Q_{2d}$ for $d\geq1$, where the vertex set $V(\frac{1}{2}Q_{2d})$ consists of elements of $\mathbb{Z}_{2}^{2d}$ of even Hamming weight and two vertices are adjacent if their Hamming distance is exactly two.
\end{enumerate}
Then $\mathcal {Q}(G)$ has no PST.

\begin{proof}
It has been showed \cite{Coutinho2014} that the spectra of these graphs are given by
 \begin{enumerate}[(i)]
\item $\text{Sp}(Q_{d})=\{d-2l:0\leq l\leq d\}$;

\item $\text{Sp}(\overline{mK_{2}})=\{2m-2, -2, 0\}$;

\item $\text{Sp}(\frac{1}{2}Q_{2d})=\left\{
\left( {\begin{array}{*{20}{c}}
{2d}\\
2
\end{array}} \right)-2l(2d-l): 0\leq l\leq d
\right\}$.
\end{enumerate}
Sine these graphs are distance regular and all eigenvalues of these graphs are integers. Hence, if $G$ is each one of these graphs, then $\mathcal {Q}(G)$ has no PST by Theorem 4.2.
\end{proof}

\section{PGST on $\mathcal{Q}$-graph}

Recall that the graph possessing PST is rarely. According to Theorem 4.2, there is no PST in $\mathcal{Q}(G)$ for any integral $r$-regular graph with $r\geq2$. Therefore, it is worthwhile for us to further investigate PGST on $\mathcal{Q}$-graph of regular graphs. In this section, we show that $\mathcal{Q}$-graph has PGST under some mild conditions. At the same time, some new families of graphs existing PGST are exhibited at the end of this section.

\paragraph{Theorem 5.1} Let $G$ be an $r$-regular connected graph of order $n$ with $m$ edges and $r\geq2$. Also let $ \lambda_{0}>\lambda_{1}>\cdots>\lambda_{d}$ be all distinct eigenvalues of $G$ and $g$ is the greatest common divisor as stated in Theorem 2.1. Assume that $G$ has PST at time $t=\frac{\pi}{g}$ between vertices $u$ and $v$.

\begin{enumerate}[(i)]
\item If $G$ is non-bipartite, then $\mathcal {Q}(G)$ has PGST between $u$ and $v$.

\item If $G$ is bipartite and $\frac{r}{g}$ is even, then $\mathcal {Q}(G)$ has PGST between the vertices $u$ and $v$.

\end{enumerate}

\begin{proof}
For the sake of convenience, we let $S=\text{supp}_{A(G)}(u)$ and $\Delta_{\lambda_i}=\sqrt{(\lambda_i+r)^2+4}$ for $i\in [d]\cup\{0\}$.

(i) $G$ is a connected non-bipartite graph. In this case, from the (i) of Proposition 3.4, we have

\begin{equation*}
\begin{split}
{{{(e_v^{m + n})}^T}\exp( - it{A_{\mathcal {Q}(G)}})e_u^{m + n}}
 &= {e^{ - it\frac{{r - 2}}{2}}}\sum\limits_{i = 0}^{d} {{e^{ - it\frac{{{\lambda_i}}}{2}}}} e_v^T{E_{{\lambda_i}}}{e_u}\left(\cos\frac{{{\Delta _{{\lambda_i}}}t}}{2} + \frac{{{\lambda_i} + r - 2}}{{{\Delta _{{\lambda_i}}}}}i\sin\frac{{{\Delta _{{\lambda_i}}}t}}{2}\right)\\
&={e^{ - it\frac{{r - 2}}{2}}}\sum\limits_{{\lambda_i} \in S} {{e^{ - it\frac{{{\lambda_i}}}{2}}}} e_v^T{E_{{\lambda_i}}}{e_u}\left(\cos\frac{{{\Delta _{{\lambda_i}}}t}}{2} + \frac{{{\lambda_i} + r - 2}}{{{\Delta _{{\lambda_i}}}}}i\sin\frac{{{\Delta _{{\lambda_i}}}t}}{2}\right),
\end{split}
\end{equation*}
where the last equality holds as $e_v^T{E_{{\lambda_i}}}{e_u}\neq0$ if and only if $\lambda_{i} \in \text{supp}_{A(G)}(u)$.

Observe that $|e^{-it\frac{r-2}{2}}|=1$ for any $t$. To prove $\mathcal {Q}(G)$ admitting PGST between the vertices $u$ and $v$, we only need to find a time $t_0$ such that
\begin{equation*}
\left|\sum\limits_{{\lambda_i} \in S}^{} {{e^{ - it_0\frac{{{\lambda_i}}}{2}}}} e_v^T{E_{{\lambda_i}}}{e_u}\left(\cos\frac{{{\Delta _{{\lambda_i}}}t_0}}{2} + \frac{{{\lambda_i} + r - 2}}{{{\Delta _{{\lambda_i}}}}}\sin\frac{{{\Delta _{{\lambda_i}}}t_0}}{2}\right)\right| \approx 1.
\end{equation*}
Since $G$ admits PST between $u$ and $v$ at time $t=\frac{\pi}{g}$. Then, it follows from Theorem $2.1$ that all the elements of $S$ are integers. Thus, $\mathcal {Q}(G)$ has no PST between $u$ and $v$, by Corollary $4.3$. Bear in mind that $G$ is connected, $|S|\geq 2$. Furthermore, assume that $r\neq \lambda_{j}\in S$ for some $j\in [d]$. Then $\lambda_{j}$ is integer and $E_{\lambda_{j}}e_{u}\neq 0$, which implies that $\lambda_{j\pm} \in \text{supp}_{A(\mathcal {Q}(G))}(u)$. Set $\Delta _{\lambda_{j}}=\sigma_{j}\sqrt{\theta_{j}}$ for each $\lambda_{j}\in S$, where $\sigma_{j}$, $\theta_{j}\in \mathbb{Z^+}$ and $\theta_{j}$ is the square-free part of $\Delta _{\lambda_{j}}^{2}$.
  Since $(\lambda_{j}+r)^{2}<(\lambda_{j}+r)^{2}+4<(\lambda_{j}+r+2)^{2}$ for $\lambda_{j}\neq-r$. Then, from Lemma 2.5,
\[
\cup \{\sqrt {{\theta_j}} :{\lambda_j} \in S\}
\]
is linearly independent over $\mathbb{Q}$. Theorem 2.4 implies that, there exist integers $\alpha$ and $p_{j}$ for each $\lambda_{j}\in S$, such that
\begin{equation}\label{17}
\alpha \sqrt {{\theta_j}}  - {p_j} \approx  - \frac{1}{{2g}}\sqrt {{\theta_j}}.
\end{equation}
Noting that $\theta_{i}=\theta_{j}$ for two distinct eigenvalues $\lambda_{i}, \lambda_{j}\in S$ if and only if $p_{i}=p_{j}$. Multiplying by $4\sigma_{j}$ in both sides of (\ref{17}), we obtain
\[
\Delta_{\lambda_{j}}\approx \frac{{4{\sigma_j}{p_j}}}{{4\alpha  + \frac{2}{g}}}.
\]
Now let $t_0=(4\alpha+\frac{2}{g})\pi$. Then one has
\[\cos \frac{{{\Delta _{{\lambda_j}}}t_0 }}{2} \approx \cos \frac{{\frac{{4{\sigma_j}{p_j}}}{{4\alpha  + \frac{2}{g}}}(4\alpha  + \frac{2}{g})\pi}}{2} = \cos 2{\sigma_j}{p_j}\pi = 1
\]
for integers ${\sigma_j}, {p_j}$.
Hence,
\begin{equation*}
\begin{split}
\left|\sum\limits_{{\lambda_j} \in S} {{e^{ - it_0 \frac{{{\lambda_j}}}{2}}}} e_v^T{E_{{\lambda_j}}}{e_u}(\cos\frac{{{\Delta _{{\lambda_j}}}t_0 }}{2} + \frac{{{\lambda_j} + r - 2}}{{{\Delta _{{\lambda_i}}}}}i\sin\frac{{{\Delta _{{\lambda_j}}}t_0 }}{2})\right|
 &\approx \left|\sum\limits_{{\lambda_j} \in S} {{e^{ - it_0 \frac{{{\lambda_j}}}{2}}}} e_v^T{E_{{\lambda_j}}}{e_u}\right|\\
& \approx \left|\sum\limits_{{\lambda_j} \in S} {{e^{ - i\frac{\pi }{g}{\lambda_j}}}} e_v^T{E_{{\lambda_j}}}{e_u}\right|=1,
\end{split}
\end{equation*}
where the last equality follows as $G$ admits PST between the vertices $u$ and $v$ at time $\frac{\pi}{g}$.\\

(ii) $G$ is a connected bipartite graph. According to the (ii) of Proposition 3.4, we only need to prove that
\begin{equation*}
\left|\sum\limits_{{\lambda_j} \in S\backslash \{-r\}} {{e^{ - it\frac{{{\lambda_j}}}{2}}}} e_v^T{E_{{\lambda_j}}}{e_u}(\cos\frac{{{\Delta _{{\lambda_j}}}t}}{2} + \frac{{{\lambda_j} + r - 2}}{{{\Delta _{{\lambda_j}}}}}i\sin\frac{{{\Delta _{{\lambda_j}}}t}}{2}) + {e^{ - it0}}e_v^T{E_{-r}}{e_u}\right| \approx 1
\end{equation*}
for some time $t$. Setting $t_0=(4\alpha+\frac{2}{g})\pi$. Similar to the discussion in (i), we obtain
\begin{equation}\label{18}
\begin{split}
&\left|\sum \limits_{{\lambda_j} \in S\backslash \{-r\}} {{e^{ - it_0 \frac{{{\lambda_j}}}{2}}}} e_v^T{E_{{\lambda_j}}}{e_u}(\cos\frac{{{\Delta _{{\lambda_j}}}t_0 }}{2} + \frac{{{\lambda_j} + r - 2}}{{{\Delta _{{\lambda_j}}}}}i\sin\frac{{{\Delta _{{\lambda_j}}}t_0 }}{2})+ {e^{ - it_0 0}}e_v^T{E_{-r}}{e_u}\right|\\
&\approx\left|\sum\limits_{{\lambda_j} \in S\backslash \{-r\}} {{e^{ - it_0 \frac{{{\lambda_j}}}{2}}}} e_v^T{E_{{\lambda_j}}}{e_u} + {e^{ - it_0 0}}e_v^T{E_{-r}}{e_u}\right| .
\end{split}
\end{equation}
Notice that ${e^{-it_00}}=e^{-it_0\frac{-r}{2}}$ as $\frac{r}{g}$ is an even integer. It follows from (\ref{18}) that
\begin{equation*}
\begin{split}
 \left|\sum\limits_{{\lambda_j} \in S\backslash \{-r\}} {{e^{ - it_0 \frac{{{\lambda_j}}}{2}}}} e_v^T{E_{{\lambda_j}}}{e_u} + {e^{ -it_0 0 }}e_v^T{E_{-r}}{e_u}\right|
 &= \left|\sum\limits_{{\lambda_j} \in S\backslash \{-r\}} {{e^{ - it_0 \frac{{{\lambda_j}}}{2}}}} e_v^T{E_{{\lambda_j}}}{e_u} + {e^{ - it_0 \frac{-r}{2}}}e_v^T{E_{-r}}{e_u}\right|\\
 &= \left|\sum\limits_{{\lambda_j} \in S} {{e^{ - it_0 \frac{{{\lambda_j}}}{2}}}} e_v^T{E_{{\lambda_j}}}{e_u}\right|\\
 & = \left|\sum\limits_{{\lambda_j} \in S} {{e^{ - i\frac{\pi }{g}{\lambda_j}}}} e_v^T{E_{{\lambda_j}}}{e_u}\right|
 = 1,
\end{split}
\end{equation*}
completing the proof this theorem.
\end{proof}

We shall exhibit some new families of graphs existing PGST at the end of this section. It is well-known that distance regular graphs have rich combinatorial structures and play a significant role in graph theory and combinatorial mathematics. In quantum walks, an interesting result is the following: the eigenvalue support of each vertex in a distance regular graph $G$ equals the set of all distinct eigenvalues of $G$ (see \cite{Coutinho2014}). It reduces the work to determine whether the eigenvalues are in the eigenvalue support of some vertex and helps us to find $g$ by a simply calculation.

\paragraph{Example 2.}
\begin{enumerate}[(i)]
  \item Note that cocktail party graph $\overline{mK_{2}}$ is a non-bipartite $(2m-2)$-regular graph with $g=2$ (see Example 1) and it admits PST between antipodal vertices $u$ and $v$ at time $\frac{\pi}{2}$ in \cite{Coutinho2014}. Hence, $\mathcal {Q}(\overline{mK_{2}})$ has PGST between $u$ and $v$;
  \item Meixner graph in \cite{Coutinho2014} is a non-bipartite distance regular graph of order 1344. Its spectrum is $\{176, 44, 8, -4 ,-16\}$ and $r=176$. By calculations, $g=12$. It was proved in \cite{Coutinho2014} that PST happens between antipodal vertices $u$ and $v$ at time $\frac{\pi}{12}$. Hence, $\mathcal {Q}$-graph of \text{Meixner graph} admits PGST between $u$ and $v$.
  \item The $d$-cube $Q_{d}$ is a bipartite distance regular graph with $r=d$ and $g=2$. It is well-known that $Q_{d}$ has PST between antipodal vertices $u$ and $v$ at time $\frac{\pi}{2}$. Then $\mathcal {Q}$-graph of $Q_{d}$ admits PGST between $u$ and $v$ whenever $d$ is multiple of four.
\end{enumerate}

\section{Concluding remarks}

This paper mainly focuses on PST and PGST in $\mathcal {Q}$-graphs of regular graphs. Studying compound graph by some graph operations is interesting and meaningful towards more larger graphs. We give spectral decomposition of adjacency matrix of $\mathcal {Q}$-graph of a regular graph. Applying these results, we discuss the existence of PST and PGST on $\mathcal {Q}$-graph of a regular graph. Finally, applying these results to some distance regular graphs, we also exhibit many new families of $\mathcal{Q}$-graphs having no PST, but admitting PGST.

The investigation in this article avoids the case of irregularity. Therefore, we put forward the following problem at the end of this article:
 \emph{For a non-regular graph, determine whether or not the corresponding $\mathcal {Q}$-graph has PST or PGST.}\\
\\
\textbf{Acknowledgements} This work was in part supported by the National Natural Science Foundation of China (No. 11801521).

\end{document}